\documentclass[a4paper,10pt]{article}

\usepackage{mathptmx}
\usepackage[english]{babel}

\usepackage{amsmath}
\usepackage{amsfonts}
\usepackage{amssymb}
\usepackage{amsthm}
\usepackage{graphicx}
\usepackage{pdfpages}

\newcommand{\setx}[1]{ \{ #1 \} }

\newcommand{\mytext}[1]{ \: \textrm{#1} \: }
\newcommand{\mysetdescr}[2]{\left\{ #1 \: \left| \: #2 \right. \right\} }
\newcommand{\darr}{{\downarrow \,}}
\newcommand{\uarr}{{\uparrow \,}}
\newcommand{\ouarr}[1]{\left( \uarr #1 \right) \setminus \setx{#1}}
\newcommand{\odarr}[1]{\left( \darr #1 \right) \setminus \setx{#1}}

\newcommand{\myN}{\mathbb{N}}
\newcommand{\myNk}[1]{\setx{ 1, \ldots , #1 }}
\newcommand{\myNkz}[1]{\setx{ 0, \ldots , #1 }}

\newcommand\urbild[1]{#1^{-1}}

\def\A{{\cal A}}
\def\B{{\cal B}}
\def\C{{\cal C}}

\def\F{{\cal F}}

\def\J{{\cal J}}
\def\K{{\cal K}}

\def\M{{\cal M}}
\def\N{{\cal N}}
\def\P{{\cal P}}
\def\V{{\cal V}}
\def\W{{\cal W}}

\newcommand{\mf}[1]{\mathfrak{ #1 }}
\newcommand{\fC}{\mf{C}}
\newcommand{\fF}{\mf{F}}
\newcommand{\fG}{\mf{G}}

\def\BP{\begin{proof}}
\def\EP{\end{proof}}

\DeclareMathOperator{\id}{id}
\DeclareMathOperator{\Conn}{Conn}
\DeclareMathOperator{\imp}{imp}

\newcommand{\CoNull}{\Conn}
\newcommand{\Co}{\CoNull_{23}}
\newcommand{\Codr}{\CoNull_3}

\newcommand{\fCd}{\fC_3^\diamond}
\newcommand{\PLU}{\P^2_{LU}}

\newcommand{\FP}{\F(P)}
\newcommand{\fFP}{\fF(P)}
\newcommand{\fFQ}{\fF(Q)}
\newcommand{\MP}{\M(P)}

\begin{document}

\theoremstyle{plain}
\theoremstyle{plain}
\newtheorem{theorem}{Theorem}
\newtheorem{definition}{Definition}
\newtheorem{corollary}{Corollary}
\newtheorem{lemma}{Lemma}
\newtheorem{proposition}{Proposition}
\newtheorem{result}{Result}

\title{\bf Crowns as retracts}
\author{\sc Frank a Campo}
\date{{\sf acampo.frank@gmail.com}}


\maketitle

\begin{abstract}
\noindent
We investigate crowns as retracts of finite posets. We define a multigraph $\fFP$ reflecting the network of so-called improper 4-crowns contained in the extremal points of $P$, and we show that $P$ contains a 4-crown as retract iff there exists a graph homomorphism of a certain type from $\fFP$ to a multigraph $\fC$ not depending on $P$. Additionally we show that $P$ contains a retract-crown with more than four points iff the poset induced by the extremal points of $P$ contains such a retract-crown. As practical result we develop and apply criteria for the systematic investigation of crowns as retracts. Most of our results are valid for infinite posets without infinite chains, too.
\newline

\end{abstract}

\section{Introduction} \label{sec_introduction}

A poset $P$ is said to have the {\em fixed point property} iff every order homomorphism $f : P \rightarrow P$ has a fixed point. A fundamental tool in the investigation of this property are {\em retractions} and {\em retracts}. In the second half of the seventies, these objects have been developed and intensively studied by a group of scientists around Rival \cite{Bjoerner_Rival_1980,Duffus_etal_1980_DPR,Duffus_Rival_1979,Duffus_Rival_1981,Duffus_etal_1980_DRS,Nowakowski_Rival_1979,Rival_1976,Rival_1980,Rival_1982}. From the later work about the fixed point property and retracts, pars pro toto the work of Rutkowski \cite{Rutkowski_1986,Rutkowski_1989,Rutkowski_Schroeder_1994}, Schr\"{o}der \cite{Rutkowski_Schroeder_1994,Schroeder_1993,Schroeder_1995,Schroeder_1996,Schroeder_2012,Schroeder_2021,Schroeder_2022_MASoC}, and Zaguia \cite{Zaguia_2008} is mentioned here. For details, the reader is referred to \cite{Schroeder_2016,Schroeder_2012}.

{\em Crowns} can be regarded as the ``natural enemies'' of the fixed point property, because a crown being a retract of a poset directly delivers a fixed point free homomorphism. Brualdi and Dias da Silva \cite{Brualdi_DdaSilva_1997} generalized this observation by showing that $P$ does not have the fixed point property iff it contains a {\em generalized crown} as retract. 

In the present paper, we investigate crowns as retracts of finite posets. Due to a result of Duffus et al.\ \cite{Duffus_etal_1980_DPR} cited at the end of Section \ref{subsec_notation}, we can without loss of generality focus on crowns contained in the sub-poset $E(P)$ induced by the extremal points of $P$.

The main part of our article deals with 4-crowns. We call a crown $\setx{a,b,v,w}$ with minimal points $a,b$ and maximal points $v,w$ an {\em improper 4-crown} iff there exists an $x \in P$ with $a, b < x < v, w$; otherwise, we call it {\em proper}. An improper 4-crown cannot be a retract of $P$, and if a proper 4-crown $C$ is a retract of $P$, all improper 4-crowns in $E(P)$ have to be mapped to $C$ in a smart way. The surprising result is that the ability to map the {\em total} poset $P$ onto $C$ in an order preserving way depends entirely on the ability to do so for the points of improper 4-crowns contained in $E(P)$.

This result implies that the points of $P$ not related to improper 4-crowns are irrelevant for $C$ being a retract of $P$ or not. At least for the non-extremal points of $P$ not being a mid-point of an improper 4-crown, this is easily understood because they can be removed by I-retractions. A detailed discussion is contained in Section \ref{subsec_Cgraph}.

Section \ref{sec_preparation} is preparatory. In Section \ref{subsec_Cgraph}, we firstly define a multigraph $\fC$ containing certain sub-posets of $C$ as vertices. Secondly, for a poset $P$, we define a multigraph $\fFP$ reflecting the network of the improper 4-crowns in $E(P)$. In Theorem \ref{theo_surjAufC}, we show that a proper 4-crown $C \subseteq E(P)$ is a retract of $P$ iff there exists a graph homomorphism from $\fFP$ to $\fC$ providing a certain separation property. In Section \ref{subsec_imagePhi}, we analyze what we can say about $P$ and $C$ if the image of such a homomorphism is a clique in $\fC$, and in Section \ref{sec_FPP_P_EP}, we present additional results about crowns as retracts. In particular, we see in Theorem \ref{theo_crown_drei} that a crown in $E(P)$ with more than four points is a retract of $P$ iff it is a retract of $E(P)$, and we see in Corollary \ref{coro_construction} that for a poset of height greater than two, the 4-crowns being a retract of it can be found by inspecting a poset of height two. In Section \ref{sec_application}, we demonstrate how our results facilitate robust criteria and approaches for the systematic check if a proper 4-crown is a retract.

With the exception of Corollary \ref{coro_EP_FPP}, all our results are valid also for infinite posets without infinite chains. Even if the maximal sets addressed in Definition \ref{def_crownBundle} do not exist, everything remains correct if we replace $\FP$ by $\F_0(P)$ in Definition \ref{def_fFP} of the multigraph $\fFP$.

\section{Preparation} \label{sec_preparation}

\subsection{Notation and basic terms} \label{subsec_notation}

In this section, we introduce our notation and recapitulate definitions of objects being in the focus of our investigation. For all other terms of order theory and graph theory, the reader is referred to standard textbooks, e.g., \cite{Diestel_2017,Schroeder_2016}. 

In what follows, $P = (X, \leq_P)$ is a finite poset. The pairs $x, y \in P$ with $x <_P y$ are called the {\em edges of $P$}. For $Y \subseteq X$, the {\em induced sub-poset} $P \vert_Y$ of $P$ is $\left( Y, {\leq_P} \cap (Y \times Y) \right)$. To simplify notation, we identify a subset $Y \subseteq X$ with the poset $P \vert_Y$ induced by it.

For $y \in X$, we define the {\em down-set} and {\em up-set induced by $y$} as
\begin{align*}
\darr_P \; y & := \mysetdescr{ x \in X }{ x \leq_P y }, \quad
\uarr_P \; y := \mysetdescr{ x \in X }{ y \leq_P x },\end{align*}
and for $x, y \in X$, the {\em interval $[x,y]_P$} is defined as $[x,y]_P := ( \uarr_P x ) \cap ( \darr_P y )$.

We define the following subsets of $X$:
\begin{align*}
L(P) & := \mytext{the set of minimal points of } P, \\
U(P) & := \mytext{the set of maximal points of } P, \\
E(P) & := L(P) \cup U(P) \mytext{ the set of extremal points of } P, \\
M(P) & := X \setminus E(P).
\end{align*}
Following our convention to identify a subset of $X$ with the sub-poset induced by it, the minimal (maximal) points of $Y \subseteq X$ are the minimal (maximal) points of $P \vert_Y $, and we write $L(Y)$ instead of $L( P \vert_Y)$ and analogously for the other sets.

\begin{figure}
\begin{center}
\includegraphics[trim = 70 715 310 75, clip]{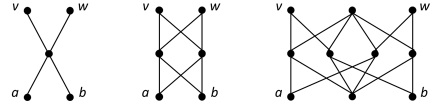}
\caption{\label{fig_Intro} Three posets $P$ of height two in which the points of $C = \setx{a,b,v,w}$ form a 4-crown in $E(P)$. In the first two posets, $C$ is an improper 4-crown, even an hourglass-crown in the leftmost one. In the rightmost poset, $C$ is a proper 4-crown, and every edge of every 4-crown in $E(P)$ belongs to an improper 4-crown.}
\end{center}
\end{figure}

For an even integer $n \geq 4$, an {\em $n$-crown} $C$ is a poset with a carrier $\setx{ c_1, \ldots , c_n }$ for which $c_1 < c_2 > c_3 < \cdots > c_{n-1} < c_n > c_1$ are the only edges of $C$. The following definition is crucial:

\begin{definition} \label{def_inner}
Let $C = \setx{a,b,v,w}$ be a 4-crown in $P$ with $\setx{a,b} = L(C)$ and $\setx{v,w} = U(C)$. We call $\J_P(C) := [a,v]_P \cap [b,w]_P$ the {\em inner of $C$}. We call $C$ a {\em proper 4-crown in $P$} iff $\J_P(C) = \emptyset$; otherwise, we call it an {\em improper 4-crown in $P$}. Furthermore, we call $C$ an {\em hourglass-crown} iff there exists an $x \in \J_P(C)$ with $\J_P(C) \subseteq ( \darr_P x ) \cup (\uarr_P x )$.
\end{definition}
Due to $[a,v]_P \cap [b,w]_P = [a,w]_P \cap [b,v]_P$, this definition is independent of the choice of the disjoint edges of $C$. As indicated by the addendum ``in $P$'', a 4-crown is proper or improper only with respect to the poset it belongs to, not as an induced sub-poset as such. However, because the ambient poset is always clear in what follows, we omit explicit reference to it in most cases. In parts of the early literature \cite{Duffus_Rival_1976,Duffus_Rival_1979,Duffus_etal_1980_DRS,Rival_1976,Rival_1982,Rutkowski_1989}, improper 4-crowns are not regarded as crowns.

Examples for the different types of 4-crowns are shown in Figure \ref{fig_Intro}. In all three posets, the points $a,b,v,w$ form a 4-crown in $E(P)$. It is an hourglass-crown in the first poset and an improper 4-crown not being an hourglass-crown in the second one. In the third poset (it is Rutkowski's \cite{Rutkowski_1989} poset P9), the points $a,b,v,w$ form a proper 4-crown. However, all {\em edges} of all crowns in $E(P)$ belong to improper 4-crowns. We will see in Theorem \ref{theo_impropCrown} that this kind of omnipresence of improper 4-crowns in $E(P)$ is typical for posets having the fixed point property.

For a homomorphism $f$ between posets or graphs, we denote by $f \vert_Z$ the pre-restriction of $f$ to a subset $Z$ of its domain. An order homomorphism $r : P \rightarrow P$ is called a {\em retraction of the poset $P$} iff $r$ is idempotent, and an induced sub-poset $R$ of $P$ is called a {\em retract of $P$} iff a retraction $r : P \rightarrow P$ exists with $r[X] = R$. For the sake of simplicity, we identify $r$ with its post-restriction and write $r : P \rightarrow R$. A poset $P$ has the fixed point property iff every retract of $P$ has the fixed point property \cite[Theorem 4.8]{Schroeder_2016}. An improper 4-crown cannot be a retract of $P$.

Finally, we want to note that according to Duffus et al.\ \cite[p.\ 232]{Duffus_etal_1980_DPR}, the retracts of $P$ in $E(P)$ contain all information about retracts of height one in $P$: Given a retraction $r : P \rightarrow R$, there exists a retraction $r' : P \rightarrow Q$ with $Q \simeq R$ and $E(Q) \subseteq E(P)$. In particular, we have $Q \subseteq E(P)$ if $R$ has height one.

\subsection{Results about homomorphisms} \label{subsec_homs}

In this technical section, we compile required tools and results about homomorphisms. Here as in all following sections, $P = (X, \leq_P)$ and $Q = (Y, \leq_Q)$ are finite posets with at least two points and without isolated points. (For the existence of homomorphisms between posets and in particular of retractions, isolated points do not play a significant role.) The sets $L(P)$, $U(P)$ and $L(Q)$, $U(Q)$ are thus disjoint antichains in $P$ and $Q$. The poset $C$ is a 4-crown with $L(C) = \setx{a,b}$ and $U(C) = \setx{v,w}$ in what follows.

Let $P$ be a poset, and let $A, B$ be a partition of $L(P)$ and $V, W$ a partition of $U(P)$. We call the mapping $f : E(P) \rightarrow C$ defined by $\urbild{f}(a) := A$, $\urbild{f}(b) := B$, $\urbild{f}(v) := V$, $\urbild{f}(w) := W$ a {\em (strict surjective) homomorphism induced by the partitions.} (In fact, the partitions allow four different induced homomorphisms, but this ambiguity will not cause problems.) 

The following proposition is the mathematical heart of this article. We need it in the proof of Lemma \ref{lemma_fourCrown} and later on in the proof of Theorem \ref{theo_crown_drei}.

\begin{proposition} \label{prop_fortsetzung}
Let $Q$ be a poset of height one and let $f : E(P) \rightarrow Q$ be a homomorphism. We define for every $x \in X$
\begin{align*}
\alpha(x) & := f \left[ L(P) \cap \darr_P x \right], \quad
\beta(x)  := f \left[ U(P) \cap \uarr_P x \right].
\end{align*}
If, for all $x \in X$ with $\alpha(x) \cap \beta(x) = \emptyset$,
\begin{align} \label{alphabeta_cond}
\begin{split}
\# \alpha(x) \geq 2 \quad & \Rightarrow \quad \# \beta(x) = 1, \\
\# \beta(x) \geq 2 \quad & \Rightarrow \quad \# \alpha(x) = 1,
\end{split}
\end{align}
then there exists a homomorphism $g : P \rightarrow Q$ with $g \vert_{E(P)} = f$.
\end{proposition}
\BP Let $x \in P$. The sets $\alpha(x)$ and $\beta(x)$ are both non-empty, each element of $\alpha(x)$ is  a lower bound of $\beta(x)$, and each element of $\beta(x)$ is an upper bound of $\alpha(x)$.

Assume $\# \alpha(x) \geq 2$. In the case of $\alpha(x) \cap \beta(x) = \emptyset$, $\beta(x)$ is a singleton according to \eqref{alphabeta_cond}. Now assume $\alpha(x) \cap \beta(x) \not= \emptyset$ and let $y \in \alpha(x) \cap \beta(x)$. Then $s \leq_Q y \leq_Q t$ for all $s \in \alpha(x)$, $t \in \beta(x)$. In the case of $y \in L(Q)$, we get $\alpha(x) = \setx{y}$ in contradiction to $\# \alpha(x) \geq 2$. Therefore, $y \in U(Q)$, thus $\beta(x) = \setx{y}$, and $\beta(x)$ is a singleton also in this case.
Together with the dual argumentation for $\# \beta(x) \geq 2$, we have seen that \eqref{alphabeta_cond} holds for {\em all} $x \in X$.

We claim that the mapping $g$ from $X$ to the carrier of $Q$ with $ g \vert_{E(P)} := f$ and 
\begin{align} \label{def_phi}
 g(x) & :=
\begin{cases}
\mytext{the single element of } \beta(x) &  \mytext{if } \# \alpha(x) \geq 2; \\
\mytext{the single element of } \alpha(x) & \mytext{if } \# \beta(x) \geq 2; \\
\mytext{the single element of } \beta(x) &  \mytext{otherwise;}
\end{cases}
\end{align}
for all $x \in M(P)$ is a homomorphism. Due to \eqref{alphabeta_cond}, the mapping is well-defined.

Let $x, y \in X$ with $x <_P y$. $ g(x) \leq_Q  g(y)$ is trivial for $x, y \in E(P)$. For $x \in M(P)$, $y \in U(P)$, we have $g(y) = f(y) \in \beta(x)$, and $g(y) = g(x)$ follows for $\setx{g(x)}= \beta(x)$. And for $\setx{g(x)} = \alpha(x)$, we have $g(x) \leq_Q g(y)$ because $g(y) \in \beta(x)$ is an upper bound of $\alpha(x)$. The case $x \in L(P)$, $y \in M(P)$ is dual.

The case $x, y \in M(P)$ remains. We have $\alpha(x) \subseteq \alpha(y)$ and $\beta(y) \subseteq \beta(x)$, and three cases are possible:
\begin{itemize}
\item $\# \alpha(x) \geq 2$: Then $\# \alpha(y) \geq 2$ , and \eqref{def_phi} delivers $\setx{g(y)} = \beta(y) \subseteq \beta(x) = \setx{g(x)}$.
\item $\# \beta(y) \geq 2$: dual to the previous case.
\item $\# \alpha(x) = 1 = \# \beta(y)$: According to \eqref{def_phi}, the equation $\# \beta(y) = 1$ yields $\beta(y) = \setx{ g(y)}$ for all cardinalities of $\alpha(y)$. In the case of $\# \beta(x) = 1$, $\beta(y) \subseteq \beta(x)$ yields $\beta(y) = \beta(x) = \setx{g(x)}$, hence $ g(x) =  g(y)$. And in the case of $\# \beta(x) \geq 2$, we have $\alpha(x) = \setx{g(x)}$, and $g(x)$ is a lower bound of $\beta(x)$, in particular also of $g(y) \in \beta(y) \subseteq \beta(x)$.
\end{itemize}
Therefore, $ g : P \rightarrow Q$ is a homomorphism with $ g \vert_{E(P)} = f$.

\EP

\begin{lemma} \label{lemma_fourCrown}
Let $C \subseteq E(P)$. There exists a retraction $f : P \rightarrow C$ iff there exists a partition $A, B$ of $L(P)$ and a partition $V, W$ of $U(P)$ fulfilling
\begin{align} \label{bedingung_Part}
\begin{split}
& \mytext{if } F \subseteq E(P) \mytext{is a four-element set with} F \cap N \not= \emptyset \mytext{ for all} \\
& \; N \in \setx{A,B,V,W}, \mytext{then} F \mytext{is not an improper 4-crown of } P
\end{split}
\end{align}
in which the points $a,b,v,w$ all belong to different sets.
\end{lemma}

\BP ``$\Rightarrow$'' holds because the existence of a retraction $f : P \rightarrow C$ implies the existence of a retraction of $P$ onto $C$ which is strict on $E(P)$.

``$\Leftarrow$'': Assume that there exists a partition $A, B$ of $L(P)$ and a partition $V, W$ of $U(P)$ with the specified properties. Let $f : E(P) \rightarrow C$ be a strict surjective homomorphism induced by these partitions.

For $x \in X$, we define the sets $\alpha(x)$ and $\beta(x)$ as in Proposition \ref{prop_fortsetzung}. Then $\alpha(x) \subseteq \setx{a,b}$ and $\beta(x) \subseteq \setx{v,w}$. Assume $\alpha(x) = \setx{a,b}$ and $\beta(x) = \setx{v,w}$. Then $x \in M(P)$, and there exist points $a' \in \urbild{f}(a) \cap \darr_P x$, $b' \in \urbild{f}(b) \cap \darr_P x$, $v' \in \urbild{f}(v) \cap \uarr_P x$, and $w' \in \urbild{f}(w) \cap \uarr_P x$. But then $F := \setx{ a', b', v', w'} \subseteq E(P)$ is an improper 4-crown of $P$ with $F \cap N \not= \emptyset$ for all $N \in \setx{A,B,V,W}$ in contradiction to \eqref{bedingung_Part}. At least one of the sets $\alpha(x)$ and $\beta(x)$ is thus a singleton, and Proposition \ref{prop_fortsetzung} yields the existence of an extension $g : P \rightarrow C$ of $f$. With an appropriate automorphism $\pi$ of $C$, the mapping $\pi \circ g$ is a retraction.

\EP

\section{Four-crowns as retracts} \label{sec_FourCrowns}

\subsection{The multigraphs $\fC$ and $\fFP$} \label{subsec_Cgraph}

\begin{figure}
\begin{center}
\includegraphics[trim = 70 630 310 74, clip]{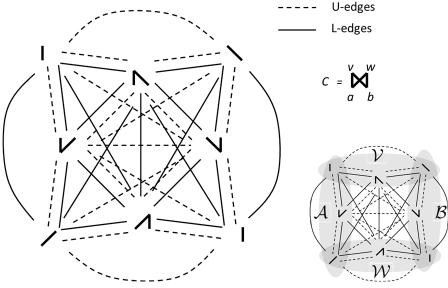}
\caption{\label{fig_CGraph} The multigraph $\fC$ with loops omitted. Explanation in text.}
\end{center}
\end{figure}

\begin{definition} \label{def_Cgraph}
With $\Co$ being the set of the connected sub-posets of $C$ with two or three points, we define $\fC$ as the undirected multigraph with vertex set $\Co$ and an edge set containing two types of edges, {\em L-edges} and {\em U-edges}, which are assigned as follows: For $S, T \in \Co$, the edge set of $\fC$ contains
\begin{itemize}
\item an  L-edge between $S$ and $T$ iff $L(C) \cap S \cap T \not= \emptyset$;
\item a U-edge between $S$ and $T$ iff $U(C) \cap S \cap T \not= \emptyset$.
\end{itemize}
\end{definition}

The multigraph $\fC$ with loops omitted is shown in Figure \ref{fig_CGraph}. The reader will be aware that the two vertical-stroke sub-posets of $C$ on top left and bottom right represent disjoint edges of the 4-crown $C$. We call the two-element sets in $\Co$ the {\em tips} of $\fC$.

According to Lemma \ref{lemma_fourCrown}, the existence of a homomorphism $f$ from $P$ onto $C$ depends entirely on the possibility to map the points of improper 4-crowns in $E(P)$ to $C$ in an appropriate way. If $f$ is strict on $E(P)$, every improper 4-crown in $E(P)$ is mapped to an element of $\Co$, and for $F$ and $G$ being improper crowns in $E(P)$ with $F \cap G \not= \emptyset$, also their images $f[F], f[G] \in \Co$ have to share points. While the graph $\fC$ is a small fixed multigraph, the complexity of the network of improper 4-crowns of $P$ is going to be reflected by their images $f[F]$ in $\fC$.

\begin{definition} \label{def_crownBundle}
We define $\C_{\imp}(P)$ as the set of improper 4-crowns contained in $E(P)$. Furthermore, we define
\begin{align*}
\MP & := \cup_{F \in \C_{\imp}(P)} \J_P(F), \\
\Xi(m) & := \left( L(P) \cap \darr m \right) \; \cup \; \left( U(P) \cap \uarr m \right) \quad \mytext{for all } m \in \MP, \\
\F_0(P) & := \mysetdescr{ \Xi(m) }{ m \in \MP }, \\
\FP & := \mytext{the set of maximal sets in } \F_0(P) \mytext{ with respect to set inclusion.}
\end{align*}
The elements of $\FP$ are called {\em 4-crown bundles}.
\end{definition}

A 4-crown bundle comprises thus the points of closely related improper 4-crowns. For an edge $x <_{E(P)} y$ in $E(P)$, the existence of a 4-crown bundle containing $x$ and $y$ is equivalent to the existence of an improper 4-crown containing these points.

\begin{definition} \label{def_fFP}
The graph $\fFP$ is the undirected multigraph with vertex set $\FP$ and an edge set containing two types of edges called {\em L-edges} and {\em U-edges} again, which are assigned as follows: For $F, G \in \FP$, the edge set of $\fFP$ contains
\begin{itemize}
\item an L-edge between $F$ and $G$ iff $L(P) \cap F \cap G \not= \emptyset$;
\item a U-edge between $F$ and $G$ iff $U(P) \cap F \cap G \not= \emptyset$.
\end{itemize}
\end{definition}

Our main tool in the investigation of 4-crowns as retracts will be a special type of homomorphisms from $\fFP$ to $\fC$ which we will introduce in Definition \ref{def_LUhom}. In order to do so, we need a covering of $\Co$:

\begin{definition} \label{def_ABVW}
We define
\begin{align*}
\A & \; := \; \mysetdescr{ S \in \Co }{ L(C) \cap S = \setx{a} }, \\
\B & \; := \; \mysetdescr{ S \in \Co }{ L(C) \cap S = \setx{b} }, \\
\V & \; := \; \mysetdescr{ S \in \Co }{ U(C) \cap S = \setx{v} }, \\
\W & \; := \; \mysetdescr{ S \in \Co }{ U(C) \cap S = \setx{w} }.
\end{align*}
\end{definition}
These four subsets of $\Co$ are indicated on the right of Figure \ref{fig_CGraph}. We note that no L-edge exists between a vertex belonging to $\A$ and a vertex belonging to $\B$. Correspondingly, vertices belonging to $\V$ and $\W$ are not connected by an U-edge.

\begin{definition} \label{def_LUhom}
Let $x, x' \in L(P)$ and $y, y' \in U(P)$ with $C' = \setx{x,x',y,y'}$ being a 4-crown. A homomorphism $\phi : \fFP \rightarrow \fC$ is called {\em $C'$-separating} iff
\begin{align} \label{phi_cond}
\mytext{for all } F \in \FP \mytext{ :} &
\begin{cases}
\; x \in F & \quad \Rightarrow \quad \phi(F) \notin \A, \\
\; x' \in F & \quad \Rightarrow \quad \phi(F) \notin \B, \\
\; y \in F & \quad \Rightarrow \quad \phi(F) \notin \V, \\
\; y' \in F & \quad \Rightarrow \quad \phi(F) \notin \W.
\end{cases}
\end{align}
\end{definition}
Let $C' = \setx{x,x',y,y'}$ be a 4-crown in $E(P)$. If there exists a $C'$-separating homomorphism $\phi$, $C'$ must be a proper 4-crown: otherwise, there exists an $F \in \FP$ with $C' \subseteq F$, and \eqref{phi_cond} yields the contradiction $\phi(F) \notin \A \cup \B \cup \V \cup \W = \Co$.

Our main result is:
\begin{theorem} \label{theo_surjAufC}
A 4-crown $C \subseteq E(P)$ is a retract of $P$ iff a $C$-separating homomorphism from $\fFP$ to $\fC$ exists.
\end{theorem}
\BP ``$\Rightarrow$'': Assume that $C$ is a retract of $P$. According to Lemma \ref{lemma_fourCrown}, there exists a partition $A, B$ of $L(P)$ and a partition $V, W$ of $U(P)$ fulfilling \eqref{bedingung_Part} in which the points of $C$ all belong to different sets. With $ g : E(P) \rightarrow C$ being the strict retraction induced by these partitions, it is easily seen that the mapping
\begin{align} \label{def_phi_von_f}
\begin{split}
\phi : \FP & \rightarrow \Co \\
F & \mapsto g[ F ],
\end{split}
\end{align}
is a homomorphism from $\fFP$ to $\fC$. We have $g(a) = a$, thus $a \in F \Rightarrow \phi(F) \notin \B$. With the corresponding results for the other points of $C$ we see that $\phi$ is $C$-separating.

``$\Leftarrow$'': Let $\phi : \fFP \rightarrow \fC$ be a $C$-separating homomorphism. We define
\begin{align*}
A_0 & \; := \; \mysetdescr{ x \in L(P)}{ \exists F \in \FP \mytext{:} x \in F \mytext{ and } \phi(F) \in \A }, \\ 
B_0 & \; := \; \mysetdescr{ x \in L(P)}{ \exists F \in \FP \mytext{:} x \in F \mytext{ and } \phi(F) \in \B }, \\ 
V_0 & \; := \; \mysetdescr{ x \in U(P)}{ \exists F \in \FP \mytext{:} x \in F \mytext{ and } \phi(F) \in \V }, \\ 
W_0 & \; := \; \mysetdescr{ x \in U(P)}{ \exists F \in \FP \mytext{:} x \in F \mytext{ and } \phi(F) \in \W }.
\end{align*}

Let $x \in A_0$ and $y \in B_0$. There exist $F, G \in \FP$ with $x \in L(P) \cap F$, $\phi(F) \in \A$, and $y \in L(P) \cap G$, $\phi(G) \in \B$. There exists no L-edge between $\phi(F)$ and $\phi(G)$ in $\fC$, and because $\phi$ is a homomorphism, we must have $x \not= y$. The sets $A_0$ and $B_0$ are thus disjoint (but one of them or both can be empty).

Because $\phi$ is $C$-separating, we cannot have $a, b \in A_0$ or $a, b  \in B_0$. Regardless of which of the points $a$ and $b$ belong to $L(P) \setminus (A_0 \cup B_0 )$ (none, one, or both), we can add them (if required) to the sets $A_0$ and $B_0$ in such a way that each of the resulting sets $A_1 \supseteq A_0$ and  $B_1 \supseteq B_0$ contains exactly one of them. Now we distribute the elements of $L(P) \setminus ( A_1 \cup B_1)$ arbitrarily on the sets $A_1$ and $B_1$. The resulting sets $A \supseteq A_1$ and $B  \supseteq B_1$ form a partition of $L(P)$ with $a$ and $b$ belonging to different sets. In the same way, we construct a partition $V, W$ of $U(P)$ with $V_0 \subseteq V$ and $W_0 \subseteq W$ in which $V$ and $W$ both contain exactly one of the points $v$ and $w$.

Let $G \in \FP$. There exists a set $\N \in \setx{\A, \B, \V, \W }$ with $\phi(G) \in \N$, say, $\N = \A$. Then $L(P) \cap G \subseteq A_0$, thus $G \cap B = \emptyset$. In the same way, we see that also for the other choices of $\N$ there exists a set $N \in \setx{A,B,V,W}$ with $G \cap N = \emptyset$, and Lemma \ref{lemma_fourCrown} delivers the existence of a retraction of $P$ onto $C$.

\EP

If $P$ is a poset of height one with $C \subseteq E(P)$, then $\FP = \emptyset$, and the trivial mapping $(\emptyset, \emptyset, \Co)$ from $\FP$ to $\Co$ is a $C$-separating homomorphism from $\fFP = (\emptyset,\emptyset)$ to $\fC$. Theorem \ref{theo_surjAufC} confirms thus the well-known result of Rival \cite{Rival_1976,Rival_1982} that a 4-crown in a poset of height one is always a retract.

The mapping $(\emptyset, \emptyset, \Co)$ is an extreme example for the fact that a separating homomorphism $\phi : \fFP \rightarrow \fC$ must not tell much about a corresponding (strict) surjective homomorphism $f : E(P) \rightarrow C$. Tracing back the construction process in the proofs of Theorem \ref{theo_surjAufC} and Lemma \ref{lemma_fourCrown}, we realize that $f$ has to send the points contained in the sets $A_0, B_0$ and $V_0, W_0$ to different points in $L(C)$ and $U(C)$, respectively. With respect to the images of the remaining points in $L(P)$ and $U(P)$, each distribution resulting in non-empty sets $A, B$ and $V, W$ is possible.

Let $Z := \left( \cup \; \FP \right) \cup \MP \cup C$ and $Q := P \vert_Z$. Because of $\fFP = \fFQ$ and Theorem \ref{theo_surjAufC}, the 4-crown $C$ is a retract of $P$ iff it is a retract of $Q$. The elements of $X \setminus Z$ are thus irrelevant for $C$ being a retract of $P$ or not. And indeed, for $x \in M(P) \setminus Z$, one of the sets $L(P) \cap \darr x $ or $U(P) \cap \uarr x$ is a singleton, and according to \cite[Ex.\ 4-5]{Schroeder_2016}, the point $x$ can be removed by a series of I-retractions. By doing so for all points in $M(P) \setminus Z$, we arrive at the poset $P' := P \vert_{E(P) \cup Z}$. The poset $P' \vert_Z$ contains $C$, the improper 4-crowns of $P$ and their inner, and the poset $P' \vert_{E(P) \setminus Z}$ has height one or is an antichain. Intuitively, it is plausible that a homomorphism mapping $P' \vert_Z$ onto $C$ can be extended to the total poset $P'$.

In the application-oriented Section \ref{sec_application}, it will be beneficial to work with a target-graph simpler than $\fC$:

\begin{corollary} \label{coro_Cdrei}
Let 
\begin{align*}
\Codr & \; := \; \mysetdescr{ S \in \Co }{ \# S = 3 }, \\
\fC_3 & \; := \; \mytext{the sub-graph of } \fC \mytext{induced by} \Codr.
\end{align*}
There exists a $C$-separating homomorphism from $\fFP$ to $\fC$ iff a $C$-separating homomorphism from $\fFP$ to $\fC_3$ exists.
\end{corollary}
\BP Only direction ``$\Rightarrow$'' has to be shown. Let $\phi$ be a $C$-separating homomorphism from $\fFP$ to $\fC$ and let $\Theta : \Co \rightarrow \Codr$ be any mapping with $\Theta(S) = S$ for all $S \in \Codr$ and $T \subset \Theta(T)$ for all two-element tips $T \in \Co \setminus \Codr$. With $\phi$, also the mapping $\Theta \circ \phi : \FP \rightarrow \Codr$ is a $C$-separating homomorphism.

\EP

\subsection{Cliques as images of homomorphisms} \label{subsec_imagePhi}

We call a non-empty subset $\K$ of $\FP$ or of $\Co$ a {\em clique} iff the vertices contained in $\K$ are all pairwise connected by an L-edge {\em and} a U-edge. The multigraph $\fC$ contains cliques with one, two, and three elements, and because each vertex in $\fC$ has an L-loop and a U-loop, every mapping from $\FP$ to a clique in $\fC$ is a homomorphism. In checking if a given poset $P$ has the 4-crown $C \subseteq E(P)$ as retract, it may thus be tempting to look if $\FP$ can be mapped to a clique in $\fC$. We see in this section that this is possible in a $C$-separating way iff an edge of $C$ does not belong to any 4-crown bundle in $\F(P)$, a result we will use in the proofs of Theorem \ref{theo_impropCrown} and Lemma \ref{lemma_LUgraph}.

\begin{lemma} \label{lemma_lemmaC}
Let $\phi : \fFP \rightarrow \fC$ be a $C$-separating homomorphism.
\begin{enumerate}
\item $C \; \not\subseteq \; \cup \; \urbild{\phi}(S)$ for every $S \in \Co$.
\item Let $S, T \in \Co$ with 
\begin{equation} \label{ab_drin}
a,b \; \in \; \left( \cup \; \urbild{\phi}(S) \right) \cap \left( \cup \; \urbild{\phi}(T) \right).
\end{equation}
Then $S = T$, if an U-edge exists between $S$ and $T$.
\end{enumerate}
\end{lemma}
\BP $C \subseteq \left( \cup \; \urbild{\phi}(S) \right)$ contradicts \eqref{phi_cond}. In the case of \eqref{ab_drin}, the vertices $S, T$ must be $\setx{a,b,v}$ or $\setx{a,b,w}$ due to \eqref{phi_cond}, and these vertices are not connected by a U-edge.

\EP

The two following propositions are formulated for 4-crown bundles, but they remain true if we exchange ``not contained in any 4-crown bundle'' against ``not contained in any improper 4-crown'' in them.

\begin{lemma} \label{lemma_fehlendeKante} 
There exists a $C$-separating homomorphism $\phi$ mapping $\FP$ to a clique $\K$ in $\fC$ iff there exists an edge of $C$ not contained in any 4-crown bundle.
\end{lemma}
\BP ``$\Rightarrow$'': We have $\# \K \leq 3$. In the case of $\K = \setx{S}$ with $S \in \Co$, Lemma \ref{lemma_lemmaC}.1 delivers $C \; \not\subseteq \; \cup \; \urbild{\phi}(S) \; = \; \cup \; \FP.
$

Assume $\K = \setx{S,T}$ with $S \not= T$, and let
\begin{align*}
Y & := \cup \; \urbild{\phi}(S), \quad \quad 
Z  := \cup \; \urbild{\phi}(T).
\end{align*}
For $C \not\subseteq Y \cup Z$, there is nothing to prove. Assume thus $C \subseteq Y \cup Z$. According to Lemma \ref{lemma_lemmaC}.2 and its dual, we have $\setx{a,b} \not\subseteq Y \cap Z$ and $\setx{v,w} \not\subseteq Y \cap Z$. Up to symmetry, two cases are possible:
\begin{itemize}
\item $a \notin Y$, $v \notin Z$: Then $a \in Z$ and $v \in Y$, and the edge $\setx{a,v}$ is not contained in any 4-crown bundle.
\item $a, v \notin Y$: Then $a, v \in Z$, hence $b \notin Z$ or $w \notin Z$ due to Lemma \ref{lemma_lemmaC}.1. In the first case, there exists no $F \in \FP$ with $b, v \in F$, and in the second case, this holds for $a, w$.
\end{itemize}

Finally, assume $\# \K = 3$ and let $T$ be the two-element tip in $\K$. As in the proof of Corollary \ref{coro_Cdrei}, we construct a $C$-separating homomorphism with image $\K \setminus \setx{T}$, and this case has already been treated.

``$\Leftarrow$'': Assume that there exists no $F \in \F(P)$ with, say, $a,v \in F$. The homomorphism
\begin{align*}
\phi(F) & :=
\begin{cases}
\setx{a,b,v}, & \mytext{if }  a \in F, \\
\setx{a,v,w}, & \mytext{otherwise.}
\end{cases}
\end{align*}
is $C$-separating because of $a \in F \Rightarrow \phi(F) \notin \A$, $v \in F \Rightarrow \phi(F) \notin \V$ and trivially $b \in F \Rightarrow \phi(F) \notin \B$ and $w \in F \Rightarrow \phi(F) \notin \W$.

\EP

We conclude
\begin{corollary} \label{coro_clique}
Assume that $\fFP$ is a complete graph. $C$ is a retract of $P$ iff there exists an edge of $C$ not contained in any 4-crown bundle.
\end{corollary}

\begin{figure}
\begin{center}
\includegraphics[trim = 70 719 290 70, clip]{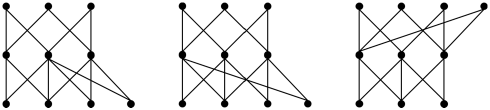}
\caption{\label{fig_P345} Rutkowski's posets P3, P4, and P5.}
\end{center}
\end{figure}

By means of this corollary we see that the posets P3-P5 of Rutkowski \cite{Rutkowski_1989} shown in Figure \ref{fig_P345} do not have a 4-crown as retract. The graph $\fFP$ of each of them is a complete graph with three vertices, but all edges of crowns contained in their extremal points belong to improper 4-crowns.

A consequence of the corollary is that a poset $P$ with at most three maximal points and at most three minimal points has a 4-crown $C \subseteq E(P)$ as retract iff one of the edges of $C$ is not contained in an improper 4-crown. With this criterion, we see immediately that also the posets P1, P2, and P6-P10 of Rutkowski \cite{Rutkowski_1989} do not have a 4-crown as retract. (P9 is shown in Figure \ref{fig_Intro} on the right.) 

\section{Additional results about crowns as retracts} \label{sec_FPP_P_EP}

Our main result about retract-crowns with more than four points is

\begin{theorem} \label{theo_crown_drei}
Let $Q$ be a poset of height one and let $f : E(P) \rightarrow Q$ be a homomorphism. There exists a homomorphism $g : P \rightarrow Q$ with $g \vert_{E(P)} = f$ if $Q$ does not contain a 4-crown. In particular, if $Q \subseteq E(P)$ is a crown with more than four points, then $Q$ is a retract-crown of $P$ iff it is a retract-crown of $E(P)$.
\end{theorem}
\BP Let $x \in P$. As in Proposition \ref{prop_fortsetzung}, we define
\begin{align*}
\alpha(x) & := f \left[ L(P) \cap \darr_P x \right], \quad
\beta(x)  := f \left[ U(P) \cap \uarr_P x \right].
\end{align*}
The absence of 4-crowns in $Q$ ensures \eqref{alphabeta_cond} for every point $x \in X$ with $\alpha(x) \cap \beta(x) = \emptyset$, and Proposition \ref{prop_fortsetzung} delivers a homomorphism $g : P \rightarrow Q$ with $g \vert_{E(P)} = f$. Direction ``$\Leftarrow$'' in the addendum is a direct consequence, and ``$\Rightarrow$'' is trivial.

\EP

For the next theorem, we need a slight variant of Rival \cite[pp.\  118-119]{Rival_1982}, \cite[Lemma 4.35]{Schroeder_2016}. The proof is a copy of the short proof of \cite[Lemma 4.35]{Schroeder_2016}:

\begin{lemma} \label{lemma_Pab}
For a poset $P$ of height one, let $C$ be an $n$-crown in $P$ and let $x, y \in C$ with $x <_P y$. $C$ is a retract of $P$ if it has the following minimality property: If $C'$ is an $n'$-crown in $P$ with $x, y \in C'$, then $n' \geq n$. 
\end{lemma}

\begin{theorem} \label{theo_impropCrown}
If $E(P)$ contains a crown with an edge not belonging to an improper 4-crown, this edge belongs to a retract-crown of $P$.
\end {theorem}

\BP Let $x <_P y$ be an edge of a crown in $E(P)$, and assume that $\setx{x, y}$ is not contained in any improper 4-crown. Under the crowns in $E(P)$ containing $x$ and $y$, select a crown $C$ with a minimal number of points. If $C$ is a 4-crown, then $C$ is a proper 4-crown, and due to Lemma \ref{lemma_fehlendeKante} and Theorem \ref{theo_surjAufC}, $C$ is a retract-crown of $P$. And if $C$ is not a 4-crown, then $C$ is according to Lemma \ref{lemma_Pab} a retract-crown of $E(P)$ with more than four points, and $C$ is a retract-crown of $P$ due to Theorem \ref{theo_crown_drei}.

\EP

We conclude that for $P$ having the fixed point property, every edge of every crown in $E(P)$ must belong to an improper 4-crown in $E(P)$. (Examples are the rightmost poset in Figure \ref{fig_Intro} and the posets in Figure \ref{fig_P345}.) For a poset $P$ of height two, this necessary condition can be sharpened: every edge of every crown in $E(P)$ must belong to an hourglass-crown, i.e., to an induced sub-poset isomorphic to the leftmost poset in Figure \ref{fig_Intro}. Otherwise, the poset $P$ contains an induced sub-poset isomorphic to the second one in Figure \ref{fig_Intro}, and $P$ does not have the fixed point property \cite[Ex. 4-7]{Schroeder_2016}.

In the following corollary and in Definition \ref{def_PLU}, we assume that every edge of every crown in $E(P)$ belongs to an improper 4-crown. According to Theorem \ref{theo_impropCrown}, this assumption is not a substantial restriction because it excludes only posets for which we already know that they contain a crown as retract.

\begin{corollary} \label{coro_construction}
Let $P$ be a poset in which every edge of every crown in $E(P)$ belongs to an improper 4-crown. There exists a poset $R$ of height two with the following properties: 
\begin{itemize}
\item The 4-crown bundles of $R$ are exactly the 4-crown bundles of $P$.
\item A 4-crown $C \subseteq E(P)$ is a retract of $P$ iff it is a retract-crown of $R$, too.
\end{itemize}
\end{corollary}
\BP We enumerate the 4-crown bundles of $P$ by $F_1, \ldots, F_n$. With $Y := \cup_{i=1}^n \; F_i$, we can assume $Y \cap \myN = \emptyset$. Now we define $R = ( Y \cup \myNk{n}, \preceq )$ by setting $x \preceq y$ iff
\begin{align*}
(x, y) & \in Y \times Y \quad \quad \quad \quad \quad \quad \quad \; \; \mytext{ and } x \leq_P y, \\
(x, y) & \in \myNk{n} \times Y \quad \quad \quad \quad \; \mytext{ and } y \in U(F_x), \\
(x, y) & \in Y \times \myNk{n} \quad \quad \quad \quad \; \mytext{ and } x \in L(F_y), \\
(x, y) & \in \myNk{n} \times \myNk{n} \quad \mytext{ and } x = y.
\end{align*}
Because of $\fF(R) = \fFP$ and Theorem \ref{theo_surjAufC}, the proposition is shown.

\EP

It are thus posets of height one and two which are in the focus of the investigation of retract-crowns. If a crown $C \subseteq E(P)$ contains more than four points, it is the poset $E(P)$ of height one which has to be inspected due to Theorem \ref{theo_crown_drei}, and if $C$ is a proper 4-crown, Corollary \ref{coro_construction} tells us that it is enough to analyze a poset of height two.

\begin{corollary}[{\cite[p.\ 275]{Duffus_Rival_1976}}] \label{coro_EP_FPP}
If $E(P)$ does not contain a crown, then $P$ is I-dismantlable to a singleton.
\end{corollary}
\BP According to Rival \cite{Rival_1976}, \cite[Theorem 4.37]{Schroeder_2016}, $E(P)$ is dismantlable to a singleton by a sequence of I-retractions, retracting in the first step, say, $\ell \in L(P)$ to $u \in U(P)$. We have $\uarr_{E(P)} \ell = \setx{\ell,u}$ due to the definition of an I-retraction, hence $(\uarr_P \ell) \cap U(P) = \setx{u}$, and $P$ is I-dismantlable to $P' := P \vert_{X \setminus \big( [\ell,u]_P \setminus \setx{u} \big)}$ according to \cite[Ex.\ 4-5]{Schroeder_2016}. Due to $E(P') = E(P) \setminus \setx{\ell}$ we can continue in processing the sequence of I-retractions of $E(P)$ in this way, and everything is shown.

\EP

\section{Application} \label{sec_application}

In this section, we demonstrate how our results facilitate robust criteria and approaches for the systematic check if a proper 4-crown is a retract.

\subsection{The poset class $\P^2_{LU}$} \label{subsec_classP2LU}

\begin{definition} \label{def_PLU}
We define $\P^2_{LU}$ as the class of finite posets $P$ of height two without isolated points for which
\begin{itemize}
\item vertices are connected in $\fFP$ by an L-edge iff they are connected by a U-edge;
\item every edge of every crown in $E(P)$ is contained in an improper 4-crown.
\end{itemize}
\end{definition}
For $P \in \P^2_{LU}$, the graph $\fFP$ can be treated as a simple graph with a single type of edges only. Moreover, no homomorphism can map adjacent 4-crown bundles to vertices in $\fC$ diagonally opposed. With Corollary \ref{coro_Cdrei}, we understand that a 4-crown $C \subseteq E(P)$ is a retract of $P \in \PLU$ iff a $C$-separating homomorphism $\phi: \fFP \rightarrow \fCd$ exists, where $\fCd$ is the 4-cycle we get by removing the L-edge $\setx{a,b,v}$-$\setx{a,b,w}$ and the U-edge $\setx{a,v,w}$-$\setx{b,v,w}$ from the multigraph $\fC_3$ defined in Corollary \ref{coro_Cdrei}. In particular, if $S, T \in \Codr$ are diagonally opposed in $\fCd$, we must have
\begin{equation} \label{nix_drin}
\left( \cup \; \urbild{\phi}(S) \right) \cap \left( \cup \; \urbild{\phi}(T) \right) \; = \; \emptyset
\end{equation}
for every $C$-separating homomorphism $\phi: \fFP \rightarrow \fCd$.

\begin{lemma} \label{lemma_LUgraph}
Let $P \in \P^2_{LU}$ and let $C \subseteq E(P)$ be a 4-crown. Every $C$-separating homomorphism $\phi$ from $\fFP$ to $\fCd$ maps the 4-crown bundles containing an edge of $C$ surjectively onto $\fCd$. Moreover, 4-crown bundles containing edges of $C$ are mapped to the same point of $\fCd$ iff they contain the same edge of $C$, and 4-crown bundles containing disjoint edges of $C$ are mapped to diagonally opposed points of $\fCd$.
\end{lemma}
\BP There exist 4-crown bundles $F_1, F_2, F_3, F_4 \in \FP$, each containing an edge of $C$ and together containing all edges of $C$. (We do {\em not} assume here that these 4-crown bundles are mutually different!). There exist $m_i \in \M(P)$ with $F_i = \Xi(m_i)$, $i \in \myNk{4}$, and with setting $Y := \cup_{i=1}^4 F_i \cup \setx{m_i}$ and $Q := P \vert_Y$, we have $\F(Q) = \setx{F_1, F_2, F_3, F_4}$ because 4-crown bundles are maximal in $\F_0(P)$. Therefore, $Q \in \PLU$ and $\fF(Q) = \fFP \vert_{\F(Q)}$, and due to Lemma \ref{lemma_fehlendeKante}, $\phi$ does not map $\F(Q)$ to a clique.

There exist thus 4-crown bundles in $\F(Q)$, say, $F_1$ and $F_4$, which are mapped to vertices $S_1$ and $S_4$ diagonally opposed in $\fCd$. Due to \eqref{nix_drin}, they contain disjoint edges of $C$ and no additional  point of $C$, and the two remaining disjoint edges of $C$ are contained in $F_2$ and $F_3$ (which now must be different because $C$ is proper!). Mapping $F_2$ or $F_3$ to $S_1$ or $S_4$ is not possible because there is no edge between $S_1$ and $S_4$, and mapping them to the same vertex of $\fCd$ contradicts Lemma \ref{lemma_lemmaC}.1. Therefore, $\phi$ maps $\F(Q)$ onto $\fCd$. The addenda follow from the way how $\F(Q)$ is mapped onto $\fCd$.

\EP

\begin{figure}
\begin{center}
\includegraphics[trim = 70 720 465 70, clip]{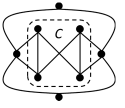}
\caption{\label{fig_PLU_forbidden} The forbidden sub-poset addressed in Lemma \ref{lemma_PLU_forbidden}.}
\end{center}
\end{figure}

This result delivers two necessary conditions for a proper 4-crown $C \subseteq E(P)$ being a retract of $P \in \PLU$ which are strict and easy to check:

\begin{lemma} \label{lemma_PLU_forbidden}
Let $P \in \P^2_{LU}$ with a proper 4-crown $C \subseteq E(P)$ and define $\C_2$ as the set of 4-crown bundles in $\FP$ containing at least two points of $C$. There exists a $C$-separating homomorphism from $\fFP \vert_{\C_2}$ to $\fCd$ iff
\begin{itemize}
\item there exists no 4-crown bundle $F \in \C_2$ with $a, b \in F$ or $v, w \in F$;
\item $P$ does not contain an induced sub-poset containing $C$ which is isomorphic to the poset in Figure \ref{fig_PLU_forbidden}.
\end{itemize}
\end{lemma}
\BP ``$\Rightarrow$: For every $F \in \C_2$, there exists a point $m_F \in \MP$ with $F = \Xi(m_F)$. With $Y := \cup_{F \in \C_2} \left( F \cup \setx{m_F}\right)$, the poset $Q := P \vert_Y$ belongs to $\PLU$ with $\fF(Q) = \fFP \vert_{\C_2}$.

Let $\psi : \fF(Q) \rightarrow \fCd$ be a $C$-separating homomorphism. According to Lemma \ref{lemma_LUgraph}, $\psi$ maps the 4-crown bundles in $\C_2$ containing an edge of $C$ surjectively onto $\fCd$. A 4-crown bundle in $\C_2$ containing $a, b$ or $v, w$ is adjacent to all of them and cannot be mapped to $\fCd$ by $\psi$. Furthermore, a structure as in the figure prevents 4-crown bundles containing disjoint edges of $C$ from being mapped to diagonally opposed vertices of $\fCd$.

``$\Leftarrow$'': Due to the first assumption, each element of $\C_2$ contains an edge of $C$ and no additional point of $C$. Therefore, the only thing we have to show is $F \cap G = \emptyset$ for all $F, G \in \C_2$ containing disjoint edges of $C$. There exist points $m, n \in \MP$ with $F = \Xi(m)$ and $G = \Xi(n)$, and we have $m \not= n$ because $C$ is proper. In the case of $F \cap G \not= \emptyset$, there exist points $z \in L(P) \cap F \cap G$ and $z' \in U(P) \cap F \cap G$. Because no element of $\C_2$ contains $\setx{a,b}$ or $\setx{v,w}$, neither $z$ nor $z'$ is contained in $C$. But then, the poset induced by the points $m, n, z, z'$ and $C$ is isomorphic to the poset in the figure.

\EP

\begin{figure}
\begin{center}
\includegraphics[trim = 70 610 190 70, clip]{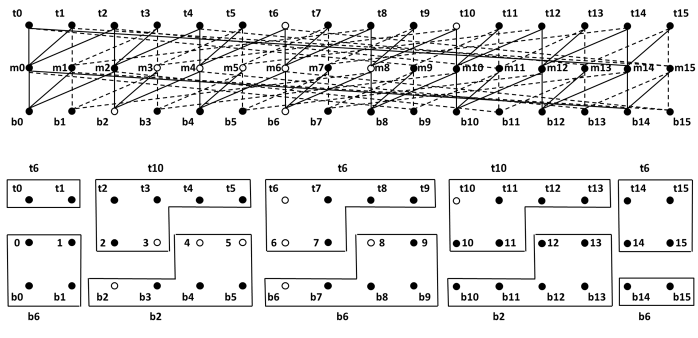}
\caption{\label{fig_X16} Top: Schr\"{o}der's poset $X_{16}$. The hollow dots $b_2, b_6, t_6, t_{10}$ form a proper 4-crown $C$, and the hollow dots $m_3, m_4, m_5, m_6, m_8$ are the mid-points of the 4-crown bundles forming the set $G_0$ used in text. For better readability, the indices of point labels are not lowered and the edges of the poset are partly dotted. In the bottom part, a retraction of $X_{16}$ to $C$ is indicated: the points within a box are mapped to the point labeling the box.}
\end{center}
\end{figure}

\subsection{Schr\"{o}der's automorphic posets}\label{subsec_Xk}

Finally, we apply our methods to certain finite {\em automorphic} posets, i.e., to finite posets having a fixed point free automorphism $\alpha$. For such a poset, there exists a proper subset $Y$ of its carrier $Z$ with $Z = \cup_{i=0}^I \alpha^i[Y]$ for an integer $I \in \myN$. We show by means of the posets $X_K$ defined by Schr\"{o}der \cite{Schroeder_2021,Schroeder_2022_MASoC} that the combination of this property with our results may yield an effective way to decide if a 4-crown is a retract. 

Let $K \geq 16$ be an even integer. The poset $X_K$ has height two and contains $3 K$ points. Its level sets are
\begin{align*}
U(X_K) & := \setx{ t_0, \ldots , t_{K-1} }, \\
M(X_K) & := \setx{ m_0, \ldots , m_{K-1} }, \\
L(X_K) & := \setx{ b_0, \ldots , b_{K-1} },
\end{align*}
and for $k \in \myNkz{K-1}$, the covering relations are
\begin{align*}
\odarr{m_k} & :=
\begin{cases}
\setx{ b_{k}, b_{k-2}, b_{k-7} }, & \mytext{if } k \mytext{ is even}, \\
\setx{ b_{k}, b_{k-1}, b_{k-2} }, & \mytext{if } k \mytext{ is odd},
\end{cases} \\
\ouarr{m_k} & :=
\begin{cases}
\setx{ t_{k}, t_{k+1}, t_{k+2} }, & \mytext{if } k \mytext{ is even}, \\
\setx{ t_{k}, t_{k+2}, t_{k+7} }, & \mytext{if } k \mytext{ is odd}.
\end{cases}
\end{align*}
Here as in what follows, index calculations are modulo $K$. The poset $X_{16}$ is shown in Figure \ref{fig_X16}. Clearly, a fixed point free automorphism $\alpha : X_K \rightarrow X_K$ is given by $\alpha(b_{k}) := b_{k+2}, \alpha(m_k) := m_{k+2}$, and $\alpha(t_{k}) := t_{k+2}$ for all $k \in \myNkz{K-1}$.

\begin{figure}
\begin{center}
\includegraphics[trim = 75 655 300 75, clip]{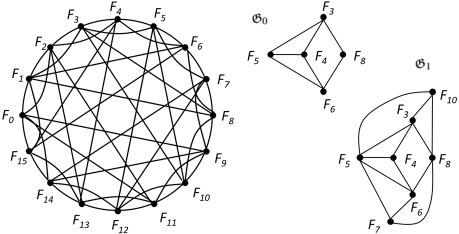}
\caption{\label{fig_ConstrX16} The graph $\fF(X_{16})$ and two subgraphs of $\fF(X_K)$ (all drawn as simple graphs with loops omitted) used in the proof that $\setx{b_2, b_6, t_6, t_{10}}$ is a retract of $X_{K}$ iff $K$ is a multiple of eight.}
\end{center}
\end{figure}

We show with our methods that for all $K \geq 16$, the poset $X_K$ has a 4-crown as retract iff $K$ is a multiple of eight (direction ``$\Leftarrow$'' has already been stated in \cite{Schroeder_2021,Schroeder_2022_MASoC}). For $K = 16$, a retraction to the 4-crown $\setx{ b_2, b_6, t_6, t_{10} }$ is indicated in the lower part of Figure \ref{fig_X16}.

Every edge in $E(X_K)$ belongs to an improper 4-crown, and by inspecting the possible combinations of even and odd, we see that the posets $X_K$ belong to $\P^2_{LU}$. For $k \in \myNkz{K-1}$, the 4-crown bundles are 
\begin{align*}
F_k & := \Xi(m_k) = 
\begin{cases}
\setx{ b_{k-7}, b_{k-2}, b_k, t_k, t_{k+1}, t_{k+2} }, & \mytext{ for even index} k, \\
\setx{ b_{k-2}, b_{k-1}, b_k, t_k, t_{k+2}, t_{k+7} }, & \mytext{ for odd index} k.
\end{cases}
\end{align*}
The graph $\fF(X_{16})$ is shown on the left of Figure \ref{fig_ConstrX16}.

Taking into account the various direct and dual symmetries of $X_K$ and the forbidden structures described in Lemma \ref{lemma_PLU_forbidden}, we see that $E(X_K)$ contains for all $K \geq 16$ two isomorphism classes of proper 4-crowns which can be retracts of $X_K$. We represent them by $\setx{ b_2, b_6, t_6, t_{10} }$ (cf.\ Figure \ref{fig_X16}) and $\setx{b_3, b_6, t_7, t_{10} }$. For $K = 16$, a third class exists with representative $\setx{ b_2, b_{10}, t_2, t_{10} }$.

We start with $C := \setx{ b_2, b_6, t_6, t_{10} }$. We define sub-graphs $\fG_k$ of $\fF(X_K)$ by setting
\begin{align*}
G_0 & := \setx{F_3, F_4, F_5, F_6, F_8}, \\
G_k & := \bigcup_{i=0}^k \alpha^i[G_0] \quad \; \; \mytext{for all } k \in \myNk{K/2}, \\
\mytext{and} \quad \fG_k & := \fF(X_K) \vert_{G_k} \quad \; \; \; \; \mytext{for all } k \in \myNkz{K/2}.
\end{align*}
The sub-graphs $\fG_0$ and $\fG_1$ are shown on the right of Figure \ref{fig_ConstrX16} with loops omitted.

The 4-crown bundles $F_3, F_4, F_6$, and $F_8$ contain the four edges of $C$, and Lemma \ref{lemma_LUgraph} tells us that a $C$-separating homomorphism $\phi_0 : \fG_0 \rightarrow \fCd$ has to map the 4-cycle $F_3$-$F_4$-$F_6$-$F_8$ surjectively onto $\fCd$. Looking at $\fG_0$ in Figure \ref{fig_ConstrX16}, we see that we must additionally have $\phi_0(F_5) = \phi_0(F_4)$.

We extend $\phi_0$ to $\fG_1$. We have $G_1 \setminus G_0 = \setx{F_7, F_{10}}$, and inspecting $\fG_1$ in Figure \ref{fig_ConstrX16}, we see that $\phi_1 \vert_{G_0} = \phi_0$ enforces $\phi_1(F_7) = \phi_0(F_6)$ and $\phi_1(F_{10}) = \phi_0(F_3)$. Setting
\begin{align*}
\rho(\phi_0(F_3)) & := \phi_0(F_4), \quad \quad
\rho(\phi_0(F_4)) := \phi_0(F_6), \\
\rho(\phi_0(F_6)) & := \phi_0(F_8), \quad \quad
\rho(\phi_0(F_8)) := \phi_0(F_3),
\end{align*}
the automorphism $\rho : \fCd \rightarrow \fCd$ is the rotation of $\fCd$ by $\pm 90^\circ$, and we have $\phi_1 \vert_{\alpha[G_0]} = \rho \circ \phi_0 \circ \urbild{\alpha}$. Proceeding in this way we see that the recursive definition
\begin{align*}
\phi_k & := 
\begin{cases}
\phi_{k-1} & \mytext{on } G_k \setminus \alpha^k[G_0], \\
\rho^k \circ \phi_0 \circ \alpha^{-k} & \mytext{on } \alpha^k[G_0],
\end{cases}
\end{align*}
is the only way to extend $\phi_0$ to a $C$-separating homomorphism $\phi_k : \fG_k \rightarrow \fCd$ for $k \in \setx{ 1, \ldots , K/2}$. However, for $\phi_{K/2}$ being well-defined, we must have $\phi_{K/2} \vert_{G_0} = \phi_0$, hence $\rho^{K/2} = \id_{\fCd}$. We conclude that $C$ is a retract of $X_K$ iff $K$ is a multiple of eight.

For $C := \setx{b_3, b_6, t_7, t_{10} }$, we start with the set $H_0 := \setx{F_3, F_5, F_6, F_7, F_8, F_{10}}$ of the 4-crown bundles containing the edges of $C$. Exactly as above, we realize that we can recursively extend a $C$-separating homomorphism $\psi_0 : \fF(X_K) \vert_{H_0}\rightarrow \fCd$ to the total graph $\fF(X_K)$ iff $K$ is a multiple of eight. Finally, in the case of $K = 16$, it are the 4-crown bundles $F_2, F_3, F_{10}, F_{11}$ which contain the four edges of the 4-crown $C := \setx{ b_2, b_{10}, t_2, t_{10} }$. They form one of the four ``diagonal stripes'' in the graph on the left of Figure \ref{fig_ConstrX16}, and by mapping the three other diagonal stripes in an appropriate way onto it, we see that there exists a $C$-separating homomorphism to $\fCd$ with pre-image sets 
\begin{equation*}
\setx{ F_0, F_2, F_4, F_{15} }, \quad
\setx{ F_1, F_3, F_5, F_6 }, \quad
\setx{ F_7, F_8, F_{10}, F_{12} }, \quad
\setx{ F_9, F_{11}, F_{13}, F_{14} }.
\end{equation*}

\noindent {\bf Acknowledgement:} I am grateful to Bernd Schr\"{o}der for his valuable comments to an early version of the manuscript and inspiring discussions.



\begin{thebibliography}{xx}

\bibitem{Bjoerner_Rival_1980} A. Bj\"{o}rner and I. Rival: A note on fixed points in semimodular lattices. {\em Discr. Math.} {\bf 29} (1980), 245--250.

\bibitem{Brualdi_DdaSilva_1997} R. A. Brualdi and J. A. Dias da Silva: A retract characterization of posets with the fixed-point property. {\em Discr. Math.} {\bf 169} (1997), 199--203.

\bibitem{Diestel_2017} R. Diestel: {\em Graph Theory.} Graduate Texts in Mathematics, Volume 173, 5\textsuperscript{th} edition, Springer (2017).

\bibitem{Duffus_etal_1980_DPR} D. Duffus, W. Poguntke, and I. Rival: Retracts and the fixed point problem for finite partially ordered sets. {\em Canad. Math. Bull.} {\bf 23} (1980), 231--236.

\bibitem{Duffus_Rival_1976} D. Duffus and I. Rival: Crowns in dismantlable partially ordered sets. In A. Hajnal and V. T. S\'{o}s  (Eds.), {\em Combinatorics}, Colloq. Math. Soc. János Bolyai, 18, Keszthely (1976), Vol. I, 271-–292.

\bibitem{Duffus_Rival_1979} D. Duffus and I. Rival: Retracts of partially ordered sets. {\em J. Austral. Math. Soc. (Series A)} {\bf 27} (1979), 495--506.

\bibitem{Duffus_Rival_1981} D. Duffus and I. Rival: A structure theory for ordered sets. {\em Discr. Math.} {\bf 35} (1981), 53--118.

\bibitem{Duffus_etal_1980_DRS} D. Duffus, I. Rival, and M. Simonovits: Spanning retracts of a partially ordered sets. {\em Discr. Math.} {\bf 32} (1980), 1--7.

\bibitem{Nowakowski_Rival_1979} R. Nowakowski and I. Rival: Fixed-edge theorem for graphs with loops. {\em J. Graph Th.} {\bf 3} (1979), 339--350.

\bibitem{Rival_1976} I. Rival: A fixed point theorem for partially ordered sets. {\em J. Comb. Th. (A)} {\bf 21} (1976), 309--318.

\bibitem{Rival_1980} I. Rival: The problem of fixed points in ordered sets. {\em Annals of Discrete Math.} {\bf 8} (1980), 283--292.

\bibitem{Rival_1982} I. Rival: The retract construction. In I. Rival (Ed.), {\em Ordered Sets. Proceedings of the NATO Advanced Study Institute held at Banff, Canada, August 28 to September 12, 1981}. Dordrecht, D.-Reidel (1982), 97--122.

\bibitem{Rutkowski_1986} A. Rutkowski: Cores, cutsets, and the fixed point property. {\em Order} {\bf 3} (1986), 257--267.

\bibitem{Rutkowski_1989} A. Rutkowski: The fixed point property for small sets. {\em Order} {\bf 6} (1989), 1--14.

\bibitem{Rutkowski_Schroeder_1994} A. Rutkowski and B. S. W. Schr\"{o}der: Retractability and the fixed point property for products. {\em Order} {\bf 11} (1994), 353--359.

\bibitem{Schroeder_2016} B. Schr\"{o}der: {\em Ordered Sets. An Introduction with Connections from Combinatorics to Topology.} Birkh\"{a}user (2016).

\bibitem{Schroeder_1993} B. S. W. Schr\"{o}der: Fixed point properties of 11-element sets. {\em Order} {\bf 2} (1993), 329--347.

\bibitem{Schroeder_1995} B. S. W. Schr\"{o}der: On retractable sets and the fixed point property. {\em Algebra Universalis} {\bf 33} (1995), 149--158.

\bibitem{Schroeder_1996} B. S. W. Schr\"{o}der: Fixed point theorems for ordered sets $P$ with $P \setminus \setx{a,c}$ as retract. {\em Order} {\bf 13} (1996), 135--146.

\bibitem{Schroeder_2012} B. S. W. Schr\"{o}der: The fixed point property for ordered sets. {\em Arab. J. Math.} {\bf 1} (2012), 529--547.

\bibitem{Schroeder_2021} B. S. W. Schr\"{o}der: The fixed point property in the set of all order relations on a finite set. {\em Order} (2021), https://doi.org/10.1007/s11083-021-09574-3.

\bibitem{Schroeder_2022_MASoC} B. S. W. Schr\"{o}der: Minimal automorphic superpositions of crowns. Arab.\ J.\ Math.\ (2022), https://doi.org/10.1007/s40065-022-00399-5.

\bibitem{Zaguia_2008} I. Zaguia: Order extensions and the fixed point property. {\em Order} {\bf 25} (2008), 267--279.

\end{thebibliography}
\end{document}